\documentclass[11pt,a4,fleqn]{article}
\usepackage{graphicx}
\usepackage{amsmath,amssymb,latexsym,graphics,epsfig}
\usepackage{hyperref}
\usepackage{color}
\usepackage{amsthm}

\newtheorem{theorem}{\bf Theorem}[section]

\newtheorem{lemma}[theorem]{\bf Lemma}

\newtheorem{emp}[theorem]{\bf Claim }

\newtheorem{conjecture}{\bf Conjecture}

\newcommand{\pf}{\noindent{\em Proof: }}
\newcommand{\epf}{\hfill\hbox{\rule{3pt}{6pt}}\\}

\numberwithin{equation}{section}

\begin{document}
\title{{\Large Monochromatic Hamiltonian Berge-cycles in colored hypergraphs}}
\author{\small  G.R. Omidi$^{\textrm{a},\textrm{b},1}$, L. Maherani$^{\textrm{a}}$\\
\small  $^{\textrm{a}}$Department of Mathematical Sciences,
Isfahan University
of Technology,\\ \small Isfahan, 84156-83111, Iran\\
\small  $^{\textrm{b}}$School of Mathematics, Institute for
Research
in Fundamental Sciences (IPM),\\
\small  P.O.Box: 19395-5746, Tehran,
Iran\\
\small \texttt{E-mails: romidi@cc.iut.ac.ir,
l.maherani@math.iut.ac.ir}}
\date {}
\maketitle \footnotetext[1] {\tt This research is partially
carried out in the IPM-Isfahan Branch and in part supported
by a grant from IPM (No. 92050217).} \vspace*{-0.5cm}

\begin{abstract} \rm{} It has been conjectured that for any fixed $r$  and sufficiently large $n$, there is a monochromatic Hamiltonian Berge-cycle in every $(r-1)$-coloring of the edges of $K_{n}^{r}$, the complete $r$-uniform hypergraph on $n$ vertices. In this paper, we show that the statement of this conjecture is true with
$r-2$ colors (instead of $r-1$ colors) by showing that there is a monochromatic Hamiltonian $t$-tight Berge-cycle in every $\lfloor\frac{r-2}{t-1}\rfloor$-edge coloring of $K_{n}^{r}$ for any fixed $r>t\geq 2$ and sufficiently large $n$.
Also, we give a proof for this conjecture when $r=4$ (the first open case). These results improve the previously known results in \cite{DGS,gyarf1,gyarf2}.


\end{abstract}

\section{\normalsize Introduction}
For given $r\geq t\geq 2$, an {\it $r$-uniform $t$-tight Berge-cycle of length $n$}, denoted by $C_n^{(r,t)}$, is an $r$-uniform hypergraph with the core sequence $v_{1}, v_{2}, \ldots ,v_{n}$ as the vertices, and distinct edges $e_{1}, e_{2}, \ldots ,e_{n}$ such that $e_{i}$ contains $v_{i}, v_{i+1},\ldots, v_{i+t-1}$ where addition is done modulo $n$. A $t$-tight Berge-cycle of length $n$ in a hypergraph with $n$ vertices is called a {\it Hamiltonian $t$-tight Berge-cycle}. This concept was introduced in \cite{DGS} to generalize Berge-cycles ($t=2$, \cite{Berge}) and tight cycles ($t=r$, \cite{hax2,RRS}). Note that, in contrast to the case $r=t=2$, for $r>t\geq 2$ a $t$-tight Berge-cycle $C_n^{(r,t)}$ is not determined uniquely and is considered as an
arbitrary choice from many possible cycles with the same triple of parameters.

Let $H$ be an arbitrary $r$-uniform hypergraph. The {\it Ramsey number} $R_{k}(H)$ is the minimum integer $n$ such that there is a monochromatic copy of $H$ in every $k$-edge coloring of $K_{n}^{r}$. The existence of such a positive integer is guaranteed by Ramsey's classical result in \cite{Ramsey}.
Recently, the Ramsey numbers of various variations of cycles in uniform hypergraphs have been studied, e.g. see \cite{hax1, hax2, OS}. Considering this problem for Berge-cycles Gy$\acute{a}$rf$\acute{a}$s et al. proposed the following conjecture:

\begin{conjecture}{\rm \cite{gyarf1}}\label{GLSS}
 Assume that $r\geq 2$ is fixed and $n$ is sufficiently large. Then every $(r-1)$-edge coloring of $K_{n}^{r}$ contains a monochromatic Hamiltonian Berge-cycle.
\end{conjecture}

This conjecture states that for a given $r\geq 2$, $R_{r-1}(C_n^{(r,2)})=n$ for sufficiently large $n$. Generalizing Conjecture \ref{GLSS} for $t$-tight Berge-cycles, Dorbec et al. proposed the following conjecture and they proved that if this conjecture is true it is best possible.

\begin{conjecture}{\rm \cite{DGS}}\label{DGC}
 Assume that $c\geq 2$, $2\leq t\leq r$, $c+t\leq r+1$  and $n$ is sufficiently large. Then every $c$-edge coloring of $K_{n}^{r}$ contains a monochromatic Hamiltonian $t$-tight Berge-cycle.
\end{conjecture}

For general cases: It is proved that the statement of Conjecture \ref{DGC} is true if we consider $ct+1\leq r$ instead of $c+t\leq r+1$ see \cite{DGS}. In \cite{gyarf1} the authors proved a weaker form of Conjecture \ref{GLSS}, which indicates that the
statement of this conjecture is true for sufficiently large $n$ with
$[\frac{r-1}{2}]$
colors instead of $r-1$
colors. In \cite{GSS} the asymptotic form of Conjecture \ref{GLSS} was proved for every $r$ using the method of Regularity Lemma. In fact, with the same assumptions the authors showed that there is a monochromatic Berge-cycle of length $(1-o(1))n$ instead of a monochromatic Hamiltonian Berge-cycle. In this paper, we improve the first two results by showing that for any fixed $r>t\geq 2$ and sufficiently large $n$, there is a monochromatic Hamiltonian $t$-tight Berge-cycle in every $\lfloor\frac{r-2}{t-1}\rfloor$-edge coloring of $K_{n}^{r}$. Clearly, this result implies that Conjecture \ref{GLSS} is true with
$r-2$ colors (instead of $r-1$ colors).

For small cases:
The case $c=2$, $t=3$ and $r=4$ of Conjecture \ref{DGC} was proved in \cite{GSS2}. In \cite{gyarf1} Conjecture \ref{GLSS} was proved for $r=3$ and an asymptotic result on this conjecture for $r=4$ was obtained using the method of Regularity Lemma.  Regarding the latter case, Gy\'{a}rf\'{a}s et al. \cite{gyarf2}, recently showed that for $n\geq 140$, in every $3$-edge coloring of $K_{n}^{4}$ there is a monochromatic  Berge-cycle  of length at least $n-10$. In the last section, we give a proof of Conjecture \ref{GLSS} for $r=4$. Our proof involves new ideas (though, it modifies certain ideas from \cite{gyarf2} at some points).

\section{\normalsize Monochromatic Hamiltonian $t$-tight Berge-cycles in colored hypergraphs}

In this section, we show that there is a monochromatic Hamiltonian $t$-tight Berge-cycle in every $\lfloor\frac{r-2}{t-1}\rfloor$-edge coloring of $K_{n}^{r}$ for any fixed $r>t\geq 2$ with $r\geq 3$ and sufficiently large $n$. This establishes the statement of Conjecture \ref{GLSS} for $r-2$ colors (instead of $r-1$ colors) and improves the former known results in \cite{DGS,gyarf1}. In order to prove our result, we need some new definitions.

Assume that $H$ is an $r$-uniform hypergraph. For a given cyclic order of $V(H)$, by a {\it consecutive $t$-vertices} we mean a subset of $V(H)$ consisting $t$ consecutive elements. The {\it shadow $t$-graph} $\Gamma_t (H)$ is a $t$-uniform hypergraph (or $t$-graph) with vertex set $V(H)$, where the edges are the sets each consisting $t$ distinct vertices for which there is an edge of $H$ containing these vertices. Let $G=\Gamma_t (H)$ and $c$ be a given $l$-edge coloring of $H$ with colors $1,2,\ldots,l$. For each edge $e=x_1x_2\ldots x_t$ of $G$, we assign a list $c(e)$ of colors of all edges of $H$ containing $x_1,x_2,\ldots, x_t$. For an edge $e$ of $G$, the color $i\in c(e)$  is called {\it $t$-good} if at least $r-t+1$ edges (of $H$) of color $i$ contain all vertices of $e$.
We consider $G$ with a new multi-coloring $c_t^{*}$ where $c_t^{*}(e)\subseteq c(e)$ is the set of all $t$-good colors for $e\in E(G)$. For $t=2$, $l=3$ and $H=K_{n}^{4}$, Gy\'{a}rf\'{a}s et al. showed that if there is a monochromatic Hamiltonian cycle $C$ in $G$ under multi-coloring $c_2^{*}$, then there is a monochromatic Hamiltonian Berge-cycle in $H$ under edge coloring $c$ (see Lemma 1 in \cite{gyarf2}). Using the same argument, we give a generalization of their result as follows:

\begin{lemma}\label{bb}
Let $r>t\geq 2$, $c$ be a given $l$-edge coloring of $H=K_{n}^{r}$ and $G=\Gamma_t (H)$. Assume that there is a monochromatic Hamiltonian tight cycle in $G$ under multi-coloring $c_t^{*}$. Then there is a monochromatic Hamiltonian $t$-tight
Berge-cycle in $H$ under $c$.
\end{lemma}

\pf Assume that $C$ is a Hamiltonian tight cycle in $G$ of color 1 (under
$c_t^{*}$) with the core sequence  $x_1,x_2,\ldots,x_n$ as the vertices. Then, following the cyclic order of vertices on $C$, suppose that $A_j$ is the set of the edges of $H$ in color 1 containing $x_j,x_{j+1},\ldots,x_{j+t-1}$. Since each $A_j$ has at least $r-t+1$ elements and no
element of $A_j$ covers more than $r-t+1$ edges  of $C$, Hall's theorem ensures the existence of
a one-to one correspondence between  all edges of $C$ (all consecutive $t$-vertices of $V(C)$) and the sets $A_j$. This clearly
defines a Hamiltonian $t$-tight Berge-cycle in $H$ under coloring $c$.
\epf

\begin{theorem}\label{Hamiltoniant-tight}
Suppose that $r>t\geq 2$ and $n\geq (r-1)\lfloor\frac{r-2}{t-1}\rfloor+2$. Then in every $\lfloor\frac{r-2}{t-1}\rfloor$-edge coloring of $K_{n}^{r}$ there is a monochromatic Hamiltonian $t$-tight Berge-cycle.
\end{theorem}

\pf Suppose to the contrary that there is no monochromatic Hamiltonian $t$-tight Berge-cycle in a given $\lfloor\frac{r-2}{t-1}\rfloor$-edge coloring of $K_{n}^{r}$. For each $1\leq i\leq \lfloor\frac{r-2}{t-1}\rfloor$, let $S_i$ be the set of all edges $e$ of $G=\Gamma_t (K_{n}^{r})$ for which $i\notin c_t^{*}(e)$. Using Lemma \ref{bb}, we may assume that the subhypergraph induced by $E(G)\setminus S_i$ in $G$ does not have a Hamiltonian tight cycle.
\begin{emp}\label{ee}
There are $(t-1)\lfloor\frac{r-2}{t-1}\rfloor+1$ vertices in $G$ so that the induced subhypergraph on these vertices in $G$ and $S_i$ have non-empty intersection, for each $1\leq i\leq \lfloor\frac{r-2}{t-1}\rfloor$.
\end{emp}
We show by induction that for each $1\leq l\leq \lfloor\frac{r-2}{t-1}\rfloor$ and any $\{S_{i_j}\}_{j=1}^{l}$ with $1\leq i_j\leq \lfloor\frac{r-2}{t-1}\rfloor$ there are $(t-1)l+1$ vertices in $G$ so that  for each $1\leq j\leq l$ the edges of the induced subhypergraph on these vertices in $G$ and $S_{i_j}$ have non-empty intersection. The case $l=1$ is trivial. Now assume that this holds for every $l<k$ where $k\leq \lfloor\frac{r-2}{t-1}\rfloor$. We verify case when $l=k$. First assume that for some $1\leq s, t \leq k$ with $s\neq t$ there are two edges $e_s\in S_{i_s}$ and $e_t\in S_{i_t}$ with $|e_s\cap e_t|\geq 2$.  By induction hypothesis there are $(t-1)(k-2)+1$ vertices in $G$ so that for each $1\leq j\leq k$ and $j\neq s,t$, the edges of the induced subhypergraph on these vertices in $G$ and $S_{i_j}$ have non-empty intersection. By adding the vertices of $e_s\cup e_t$ to these $(t-1)(k-2)+1$ vertices we get at most $(t-1)k+1$ vertices with the desired property. So we may assume that $|e_s\cap e_t|\leq 1$ for any  $1\leq s, t \leq k$ with $s\neq t$ and any two edges $e_s\in S_{i_s}$ and $e_t\in S_{i_t}$. For each $1\leq i\leq \lfloor\frac{r-2}{t-1}\rfloor$, let $G_i$ be the subhypergraph of $G$ induced by $S_i$ and $T_i$ be the set of all isolated vertices of $G_i$. Assume that $H_i$ is the subhypergraph of $G_i$ induced by $V(G_i)\setminus T_i$. We show that $\chi(H_i)>|T_i|$ for each $i\in\{i_1,\ldots,i_k\}$. Assume to the contrary that for some $i\in\{i_1,\ldots,i_k\}$ we have $\chi(H_i)\leq|T_i|$ and $C_1,C_2,\ldots,C_{\chi(H_i)}$ are the color classes of $H_i$. Let $T_i=\{t_1,t_2,\ldots,t_{|T_i|}\}$ and  $V(C_j)=\{x_{j1},x_{j2},\ldots,x_{jl_j}\}$ for $1\leq j\leq \chi(H_i)$. Consider a cyclic order of vertices of $G$ as follows:
$$ S=\{t_1,x_{11},\ldots,x_{1l_1},t_2,x_{21},\ldots,x_{2l_2},\ldots,x_{\chi(H_i)l_{\chi(H_i)}},t_{\chi(H_i)+1},t_{\chi(H_i)+2},\ldots, t_{|T_i|}\}.$$
Clearly each edge of $G$ containing an element of $T_i$ (also, each $t$-subset of any color class of $H_i$) is in $E(G)\setminus S_i$. Therefore, the set of all consecutive $t$-vertices in the cyclic order of vertices on $S$ makes a Hamiltonian tight
cycle for the subhypergraph induced by $E(G)\setminus S_i$ in $G$, a contradiction. So $\chi(H_i)>|T_i|$ for every $i\in\{i_1,\ldots,i_k\}$. Clearly for each $1\leq j\leq k$ and for any two color classes $C_s$ and $C_t$ with $s<t$ of $H_{i_j}$, there is an edge $e_{jst}\subseteq C_s\cup C_t$
in $S_{i_j}$. For such an edge $e_{jst}$ assume that $A_{jst}=(e_{jst}\cap C_s)\times (e_{jst}\cap C_t)$. By the previous argument,
$A_{jst}\cap A_{j's't'}\neq\emptyset$ if and only if $j=j',s=s'$ and $t=t'$. On the other hand, $|A_{jst}|\geq t-1$.
Therefore,
$$\sum_{j=1}^{k}{ |T_{i_j}| \choose 2}<(t-1)\sum_{j=1}^{k}{ |T_{i_j}|+1 \choose 2}\leq\sum_{j=1}^{k}\sum_{1\leq s<t\leq \chi(H_{i_j}) }|A_{jst}|,$$
which means that there is an element $(u,v)\in A_{qst}$ for some $q,s,t$ so that $\{u,v\}\nsubseteq T_{i_j}$ for each $1\leq j\leq k$. Hence, for every $j\neq q$, there is an edge $e_{i_j}$ in $S_{i_j}$ containing at least one of $u$ and $v$ as a vertex. Therefore for each $1\leq p\leq k$, the induced subhypergraph on $W=e_{qst}\cup \bigcup_{j\neq q}e_{i_j}$ and $S_{i_p}$ have non-empty intersection. Clearly, $|W|\leq (t-1)k+1$  which completes the proof of our claim.

Now, for every $1\leq i\leq \lfloor\frac{r-2}{t-1}\rfloor$ let $E_i$ be the set of all edges of color $i$ in $K_{n}^{r}$ containing all $(t-1)\lfloor\frac{r-2}{t-1}\rfloor+1$ vertices disrupted in Claim \ref{ee}. Clearly, $$\sum_{i=1}^{\lfloor\frac{r-2}{t-1}\rfloor} |E_i|\geq n-(t-1)\lfloor\frac{r-2}{t-1}\rfloor-1\geq (r-t)\lfloor\frac{r-2}{t-1}\rfloor+1.$$ On the other hand, for each $1\leq i\leq \lfloor\frac{r-2}{t-1}\rfloor$ all edges in $E_i$ contain the element $e_i\in S_i$ as a subset and $i\notin c_t^{*}(e_i)$. Hence, $|E_i|\leq r-t$ and so $\sum_{i=1}^{\lfloor\frac{r-2}{t-1}\rfloor} |E_i|\leq (r-t)\lfloor\frac{r-2}{t-1}\rfloor$, a contradiction.
\epf

The following interesting result on Conjecture \ref{GLSS} is an immediate consequence of Theorem \ref{Hamiltoniant-tight} for $t=2$.

\begin{theorem}\label{Hamiltonian}
Suppose that $n\geq r^2-3r+4$. Then in every $(r-2)$-edge coloring of $K_{n}^{r}$ there is a monochromatic Hamiltonian Berge-cycle.
\end{theorem}

\section{\normalsize Monochromatic Hamiltonian Berge-cycles in colored complete $4$-graphs}

Regarding  the case $r=4$ of Conjecture \ref{GLSS}, an asymptotic result has been obtained using the method of Regularity Lemma; see \cite{gyarf1}. Also, Gy\'{a}rf\'{a}s et al. \cite{gyarf2}, recently showed that for $n\geq 140$, in every $3$-edge coloring of $K_{n}^{4}$ there is a monochromatic  Berge-cycle of length at least $n-10$. Here, we give a proof of Conjecture \ref{GLSS} for $r=4$.

\begin{lemma}\label{a}
Suppose that $n\geq 85$ and the edges of $H=K_{n}^{4}$ are colored with three colors $1,2,3$. If there exists a vertex $v\in V(H)$ such that for some $i\in\{1,2,3\}$, at most one edge of color $i$ contains $v$, then there is a monochromatic Hamiltonian Berge-cycle in $H$.
\end{lemma}

\pf
Assume that $c$ is a $3$-edge coloring of $H$ where all edges containing $v_1=v$ are colored with colors $2$ and $3$ except possibly the edge $e_{v_1}=\{v_1,v_2,v_3,v_4\}$. Without any loss of generality, we may assume that $c(e_{v_1})\neq 3$. Consider $c'$ as the new edge coloring of $H$ such that $c'(e_{v_1})=2$ and $c'(e)=c(e)$ for any $e\in E(H)\setminus \{e_{v_1}\}$.
A new $2$-edge coloring for the $3$-uniform complete hypergraph $K$ with $n-1$ vertices $V(H)\setminus \{v_1\}$ is induced by $c'$  as follows: The edge $\{x,y,z\}$ is of color $2$ (resp. $3$) in $K$ if and only if the edge $\{v_1,x,y,z\}$ is of color $2$ (resp. $3$) under $c'$  in $H$. By Theorem 1.2 in \cite{gyarf1} there exists a monochromatic Hamiltonian Berge-cycle in $K$, say $C$. Let $x_{1},x_{2},\ldots, x_{n-1}$ be the core sequence of $C$. We consider the following cases:

\bigskip
\noindent \textbf{Case 1. } $C$ is of color $3$.\\

Without any loss of generality, we may assume that $x_1\notin\{v_2,v_3,v_4\}$. If for some non-consecutive vertices $x_{k}$ and $x_{k'}$ in $V(C)\setminus \{x_{n-1}, x_{2}\}$ the edge $\{x_1 , x_{k}, x_{k'}\}$ is of color $3$, then the cyclic order $x_{1},v_1,x_{2},\ldots, x_{n-1}$ represents a core sequence of a Hamiltonian  Berge-cycle of color $3$ in $H$. It suffices to add $v_1$ to the edge $\{x_1, x_{k}, x_{k'}\}$ and all edges of $C$ to get the edges  of a Hamiltonian  Berge-cycle of color $3$ in $H$.

If $\{ x_{n-1}, x_{1}, x_{l}\}$ (resp. $\{x_1, x_{2}, x_{l}\}$) is of color $3$, for at least two numbers $l\neq n-2, 2$ (resp. $l\neq n-1, 3$), then we have the cyclic order $v_1,x_{1},\ldots,x_{n-1}$ (resp. $x_{1},v_1,x_{2},\ldots,x_{n-1}$) representing a core sequence of a Hamiltonian Berge-cycle of color $3$ in $H$. It is sufficient to add $v_1$ to the edge $\{x_{n-1}, x_{1}, x_{l}\}\notin E(C)$ (resp. $\{x_1, x_{2}, x_{l}\}\notin E(C)$) and all edges of $C$ to get the edges of a Hamiltonian  Berge-cycle of color $3$ in $H$.

Now, we may assume that $\{x_{n-1}, x_{1}, x_{l}\}$ (resp. $\{x_1, x_{2}, x_{l}\}$) is of color $2$, for at least $n-6$ numbers $l\neq n-2, 2$ (resp. $l\neq n-1, 3$). Also, for any two non-consecutive vertices $x_{k}$ and $x_{k'}$ in $V(C)\setminus \{x_{n-1}, x_{2}\}$, the edge $\{x_1 , x_{k}, x_{k'}\}$ is of color $2$. Consider a new cyclic order $y_{1}=v_1, y_{2},\ldots ,y_{n-1},y_{n}=x_1$ for $V(H)$ such that for each $2\leq i\leq n-1$, two vertices $y_{i}$ and $y_{i+1}$ don't appear as consecutive vertices in $V(C)$ and for any $2\leq i\leq n-2$, the edge $\{y_{i},y_{i+1},x_1\}$ is of color $2$. This is possible if we set $y_{3}=x_{n-1}, y_{6}=x_{2}$ and we choose $y_2$ and $y_4$ (also $y_5$ and $y_7$) as two non-consecutive vertices in $V(C)\setminus \{x_{n-2},x_{n-1},x_{1},x_{2},x_{3}\}$ such that $\{y_{3},y_{i},x_{1}\}$ for $i=2,4$ and $\{y_{6},y_{i},x_{1}\}$ for $i=5,7$ are of color $2$. The cyclic order $y_{1}, y_{2},\ldots, y_{n}$ defines a Hamiltonian Berge-cycle of color $2$ in $H$ with the following edge assignments. Set $e_{i}=\{v_1,y_{i},y_{i+1},x_1 \}$ for $2\leq i \leq n-2$, $e_{n-1}=\{v_1,y_{p},y_{n-1},x_{1}\}$,  $e_{n}=\{v_1,y_{h},y_{k},x_1 \}$ and $e_{1}=\{v_1,y_{2},y_{l},x_1 \}$, where $y_p,y_h,y_k$ and $y_l$ are pairwise non-consecutive vertices in $V(C)\setminus \{y_{n-2},y_{n-1},y_{n},y_{1},y_{2},y_{3},y_{6}\}$ and $y_{n-1}$ and $y_{p}$ (also,  $y_2$ and $y_{l}$) are non-consecutive vertices in $V(C)$.

\bigskip
\noindent \textbf{Case 2. } $C$ is of color $2$.\\

If $\{v_2,v_3,v_4\}\notin E(C)$, then by an argument similar to that in case 1 we can see that there is a monochromatic Hamiltonian Berge-cycle in $H$.
Now, suppose that the edge $e_{1}=\{v_2,v_3,v_4\}$ appears in $E(C)$ to cover the consecutive vertices $v_2$ and $v_3$. We may assume that $x_{1}=v_2, x_{2}=v_3, x_{3},\ldots,x_{n-1}$ is the core sequence of the cycle $C$ where for each $1\leq i\leq n-1$, $e_{i}\in E(C)$ is the edge containing $x_{i}$ and $x_{i+1}$. If there are two distinct edges
$\{v_2, x_k, x_{k'}\}$ and $\{v_3, x_l, x_{l'}\}$ of color $2$ in $E(K)\setminus E(C)$,  then we consider the cyclic order of vertices of $H$ as $y_{1}=v_2, y_{2}=v_1, y_{3}=v_3, y_{4}=x_{3},\ldots,y_{n}=x_{n-1}$. The edges $f_{1}=\{v_2,v_1,x_k,x_{k'}\}$, $f_{2}=\{v_1,v_3,x_l,x_{l'}\}$ and for $3\leq i\leq n$,  $f_{i}=e_{i-1}\cup \{v_1\}$ define a Hamiltonian  Berge-cycle of color $2$ in $H$. So we may assume that for at least one of the vertices $v_2$ and $v_3$, say $v_2$, all the edges $\{v_2,x_k,x_{k'}\}\neq e_{1}, e_{n-1}$ are of color $3$ where $x_k$ and $x_{k'}$ are non-consecutive vertices of $C$. Now, we consider a new cyclic order $y_{1}=v_1, y_{2},y_{3},\ldots,y_{n-1},y_{n}=v_2$ of the vertices $V(H)$, where for any $2\leq i\leq n-1$, $y_{i},y_{i+1}$  are not consecutive vertices in $V(C)$ and for any $2\leq i\leq n-2$ the edge $\{y_{i},y_{i+1},v_2 \}$ is of color $3$. Clearly $v_{3}$ and $v_{4}$ are not consecutive vertices of the mentioned cyclic order. The following edge assignments for this cyclic order represent a Hamiltonian Berge-cycle of color $3$ in $H$, which completes the proof.
Set $f_{i}=\{v_1,y_{i},y_{i+1},v_2\}$ for $2\leq i\leq n-2$, $f_{n-1}=\{v_1,y_{p},y_{n-1},v_2 \}$, $f_{n}=\{v_1,y_{h},y_{k},v_2 \}$ and $f_{1}=\{v_1,y_{2},y_{l},v_2 \}$, where $4\leq p,h,k,l\leq n-3$, $y_{p}, y_{h}, y_{k}$ and $y_{l}$ are non-consecutive vertices in $V(C)\setminus (\{v_3,v_4\}\cup e_{n-1})$ and $y_{n-1}$ and $y_{p}$ (also,  $y_{2}$ and $y_{l}$) are non-consecutive vertices in $V(C)$.
\epf

\begin{theorem}\label{aa}
Any $3$-edge coloring of $K_{n}^{4}$ with $n\geq 85$ contains a monochromatic Hamiltonian Berge-cycle.
\end{theorem}

\pf Assume that $c$ is a $3$-edge coloring of $H=K_{n}^{4}$ with colors $1,2,3$. In \cite{gyarf2} under the same assumptions Gy\'{a}rf\'{a}s et al.  showed that if $\vert c_2^{*}(e)\vert =1$ for an edge $e$ of $G=\Gamma_2 (H)$, then there is a monochromatic Hamiltonian Berge-cycle in $H$. So suppose that  for any edge $e$ of $G$, we have $|c_2^{*}(e)|\geq 2$.

Let $v$ be an arbitrary vertex. Define $U_{12} (v)$, $U_{13} (v)$, $U_{23} (v)$ and $U_{123} (v)$ as the sets to which $v$ is connected (in the multi-coloring $c_2^{*}$) in color sets $12$, $13$, $23$ and $123$, respectively.  For $i,j,k \in \{1,2,3\}$ in some order, define
$$B_{i}=\{v \in V(G) \vert U_{ij}(v)=U_{ik}(v)=\emptyset, U_{jk}(v)\neq \emptyset \},
B_{4}=\{v \in V(G) \vert \vert U_{123}(v)\vert \geq \frac{n}{2}\}.$$
It is easy to see that for $i\neq 4$, $B_{i}$'s are pairwise disjoint and for an edge $e$ of $G$ from $B_{i}$ to $B_{j}$ where $i,j \in \{1,2,3\}$ and $i\neq j$, we have $c_2^{*}(e)=\{1,2,3\}$. Also, in \cite{gyarf2} it has been shown that if $V(G)=\bigcup _{i=1}^{4} B_{i}$, then there is a monochromatic Hamiltonian cycle for $G$ under the multi-coloring $c_2^{*}$ and so by Lemma \ref{bb} for $r=4$ and $t=2$, we conclude that there is a monochromatic Hamiltonian Berge-cycle in $H$. So suppose that $\bigcup _{i=1}^{4} B_{i}\neq V(G)$. For every $v \in V(G)\setminus \bigcup _{i=1}^{4} B_{i}$, consider $\pi(v)=\min\{\vert U_{23}(v)\vert, \vert U_{12}(v)\vert, \vert U_{13}(v)\vert \}$. We choose a vertex $x \in V(G)\setminus \bigcup _{i=1}^{4} B_{i}$ with minimum $\vert U_{123}(x)\vert$, among those with minimum $\pi(x)$. In the sequel, for simplicity we denote $U_{ij}(x)$ and $U_{123}(x)$ ($i,j \in \{1,2,3\}$) by $U_{ij}$ and $U_{123}$, respectively. Let $U= V(G)\setminus (\{x\}\cup U_{123})$ and without any loss of generality, assume that  $\vert U_{23}\vert \leq  \vert U_{12}\vert \leq \vert U_{13}\vert $. One can easily see that $U_{12} \neq \emptyset$ and $\vert U\vert \geq \lfloor \frac{n}{2}\rfloor$.

In \cite{gyarf2} it has been shown that $\vert U_{23}\vert \leq 1$. Now we show  that if $\vert U_{23}\vert=1$, then $\vert   U_{12}\vert \leq 2$. Since $\vert U\vert \geq \lfloor \frac{n}{2}\rfloor$, there are at least nine vertices in $U_{13}$. Let $S$ be a subset of $U_{13}$ of cardinality $9$.  Suppose that $u\in U_{23}$ and $T\subseteq U_{12}$ with $|T|=3$. There are twenty seven edges in $H$ each consisting $x, u$, one of the vertices in $T$ and one member of $S$. On the other hand, at most two of these edges are of color 1 (each edge has $u$ as a vertex),
at most six of them are of color 3 (each edge has exactly one of the vertices in $T$) and at most eighteen of them are of color 2 (each edge has exactly one of the vertices in $S$), a contradiction.
%In the sequel, we define a graph $\G$ with $V(\G)=V(H)$ in such a way $\G$  has a cycle $C$ of length at least $n-1$ as a subgraph. Then we show that such a cycle $C$ can be extended to a monochromatic Hamiltonian Berge-cycle in $K$.

In the sequel, we assume that $y\in U_{12}$ and $z\in U_{13}$ are fixed vertices and we define a Hamiltonian graph $\Gamma$ with $V(\Gamma)=V(H)$, in such a way that every Hamiltonian cycle $C$ of $\Gamma$  can be extended to a monochromatic Hamiltonian Berge-cycle in $H$. For this, we consider the following two cases:

\bigskip
\noindent \textbf{Case 1. }
Let $\vert U_{23}\vert =\emptyset$.\\

Let $U_{123}$ and $U_{12}$ be partitioned into  $A,B$ and $A',B'$ respectively, where $\vert B\vert \leq \vert A| \leq \vert B\vert +1$ and $\vert B'\vert \leq \vert A'| \leq \vert B'\vert +1$.  Suppose that $y\in A'$.
Consider a graph $\Gamma$ with the vertex set $V(\Gamma)=V(H)$ and the edge set $E(\Gamma)=\bigcup _{i=1}^{7}E_{i}$, where $E_{i}$s are defined as follows: 
\begin{itemize}
\item[$\bullet$]
$E_{1}=\{ uv\vert u\in U_{12}\setminus \{y\}, v\neq y, c(\{ x,z,u,v\})=1\}$.
\item[$\bullet$]
$E_{2}=\{ uv\vert u\in U_{13}\setminus \{z\}, c(\{ x,y,u,v\})=1\}$.
\item[$\bullet$]
$E_{3}=\{ zv\vert v\in V(\Gamma)\setminus A\cup A', c(\{ x,y,z,v\})=1\}$.
\item[$\bullet$]
$E_{4}=\{ yv\vert v\in A\cup A', c(\{ x,y,z,v\})=1\}$.
\item[$\bullet$]
$E_{5}=\{ yv\vert v\in U_{13}\setminus \{z\}\}$.
\item[$\bullet$]
Assume that $U_{123}=\{w_{1}, w_{2},\ldots ,w_{m}\}$ and $d_{\Gamma'}(w_{1})\leq d_{\Gamma'}(w_{2})\leq \cdots \leq d_{{\Gamma'}}(w_{m})$, where $\Gamma'$ is the graph induced by
 $\bigcup _{i=1}^{5} E_{i}$.  For $i=1,2$, assume that $e_{w_{i}v}=\{x,z,w_{i},v\}$ when $w_{i}v \in E_{1}\cup E_{4}$ and $e_{w_{i}v}=\{x,y,w_{i},v\}$ when $w_{i}v \in E_{2}\cup E_{3}$.  Since $1\in c_2^{*}(xw_{i})$ there are $r_{i}=$ max $\{3-d_{{\Gamma'}}(w_{i}),0\}$ edges $W_{i}=\{g_{i1}, \ldots , g_{ir_{i}}\} \subseteq E(H)\setminus \{e_{w_{i}v} |w_{i}v\in \bigcup_{i=1}^{4} E_{i}\}$ of color $1$ containing $x$ and $w_{i}$ for $i=1,2$. Consider the following cases:
 \begin{itemize}
 \item[i.]
  $r_1\leq 2 $. Set $E'_{6}=D_{w_{1}}=\emptyset$.  If $r_2\leq 1$, then set $ E''_{6}=D_{w_{2}}=D'_{w_{2}}=\emptyset$. Now assume that $r_{2}=2$. Let $W_{2}'=\{g_{21}, g_{22},e_{w_{2}v}\}$ where $w_{2}v\in \bigcup_{i=1}^{4} E_{i}$. Set $D'_{w_{2}}=\emptyset$, $E''_{6}=\{w_{2}t_{1}\}$ and  $D_{w_{2}}=g_{21}\setminus \{x,w_{2},t_{1}\}$ where $t_{1}\in g_{21}\setminus \{x,w_{2},v\}$. Now let $r_{2}=3$. Set $E''_{6}=\{w_{2}t_{1},w_{2}t_{2}\}$,  $D_{w_{2}}=g_{21}\setminus \{x,w_{2},t_{1}\}$ and $D'_{w_{2}}=g_{22}\setminus \{x,w_{2},t_{2}\}$ where $t_{1}\in g_{21}\setminus \{x,w_{2}\}$ and $t_{2}\in g_{22}\setminus \{x,w_{2}, t_{1}\}$ are the vertices  with maximum repetitions in $g_{21}$ and $g_{22}$.
 \item[ii.]
 $r_{1}=3$. So we have $W_{1}=\{g_{11}, g_{12},g_{13}\}$. If $r_{2}\leq 1$, then set $E''_{6}=D_{w_{2}}=D'_{w_{2}}=\emptyset$, $E'_{6}=\{w_{1}u\}$ and $D_{w_{1}}=g_{11}\setminus \{x,w_{1},u\}$ where $u \in g_{11}\setminus \{x,w_{1}\}$.
  If $r_{2}=2$, then $W_{2}=\{g_{21}, g_{22}\}$. Let $W_{2}'=\{g_{21}, g_{22},e_{w_{2}v}\}$ where $w_{2}v\in \bigcup_{i=1}^{4} E_{i}$. We may assume that $g_{13}\notin W_{2}'$. Set  $E'_{6}=\{w_{1}u\}$, $E''_{6}=\{w_{2}t_{1}\}$, $D_{w_{1}}=g_{13}\setminus \{x,w_{1},u\}$, $D_{w_{2}}=g_{21}\setminus \{x,w_{2},t_{1}\}$ and $D'_{w_{2}}=\emptyset$ so that $u \in g_{13}\setminus \{x,w_{1}\}$ and $t_{1}\in g_{21}\setminus \{x,w_{2},v\}$. Now let $r_{2}=3$. If $W_{1}\cap W_{2}= \emptyset$, then set $E'_{6}=\{w_{1}u\}$, $E''_{6}=\{w_{2}t_{1},w_{2}t_{2}\}$, $D_{w_{1}}=g_{11}\setminus \{x,w_{1},u\}$, $D_{w_{2}}=g_{21}\setminus \{x,w_{2},t_{1}\}$ and $D'_{w_{2}}=g_{22}\setminus \{x,w_{2},t_{2}\}$ so that $u \in g_{11}\setminus \{x,w_{1}\}$, $t_{1}\in g_{21}\setminus \{x,w_{2}\}$, $t_{2}\in g_{22}\setminus \{x,w_{2}, t_{1}\}$ and $t_{1}$ and $t_{2}$ have maximum repetitions in $g_{21}$ and $g_{22}$. Otherwise,  we may assume that $|g_{1i}\cap \{w_{2}\}|\geq |g_{1j}\cap \{w_{2}\}|$ for $i<j$. Choose $t_{1}\in g_{22}\setminus \{x,w_{1}, w_{2}\}$ and set  $E'_{6}=\{w_{1}w_{2}\}$, $E''_{6}=\{w_{2}w_{1},w_{2}t_{1}\}$, $D_{w_{1}}=g_{11}\setminus \{x,w_{1},w_{2}\}$, $D_{w_{2}}=g_{22}\setminus \{x,w_{1},w_{2}, t_{1}\}$ and $D'_{w_{2}}=\emptyset$. 
 \end{itemize}
 In all cases set $E_{6}=E'_{6}\cup E''_{6}$ and $D=D_{w_{1}}\cup D_{w_{2}}\cup D'_{w_{2}}$.
% 
%   If $r_1\leq 2 $ and $r_2\leq 1$, then set $E'_{6}= E''_{6}=D_{w_{1}}=D_{w_{2}}=D'_{w_{2}}=\emptyset$. If $r_1=3$ and $r_2\leq 1$, then set $E''_{6}=D_{w_{2}}=D'_{w_{2}}=\emptyset$, $E'_{6}=\{w_{1}t\}$, $D_{w_{1}}=g_{11}\setminus \{x,w_{1},t\}$, where $t\in g_{11}\setminus \{x,w_{1}\}$.
%   Now let $r_1\leq 2$, $r_2\geq 2$ and $|g_{2i}\cap \{w_{1}\}|\geq |g_{2j}\cap \{w_{1}\}|$ for $i<j$. Choose $t_{1}\in g_{21}\setminus \{x,w_{2}\}$ in such a way that $t_{1}=w_{1}$ if $w_{1}\in g_{21}$ and $t_{2}\in g_{22}\setminus \{x,w_{2}, t_{1}\}$. Now set $E'_{6}=D_{w_{1}}=\emptyset$, $E''_{6}=\{w_{2}t_{1},w_{2}t_{2}\}$, $D_{w_{2}}=g_{21}\setminus \{x,w_{2},t_{1}\}$ and $D'_{w_{2}}=g_{22}\setminus \{x,w_{2},t_{2}\}$.
%   Finally let $r_1=3$ and $r_2\geq 2$.
%    First let $W_{1}\cap W_{2}=\emptyset$. Choose $t\in g_{11}\setminus \{x,w_{1}\}$, $t_{1}\in g_{21}\setminus \{x,w_{2}\}$, $t_{2}\in g_{22}\setminus \{x,w_{2}, t_{1}\}$ and set $E'_{6}=\{w_{1}t\}$, $E''_{6}=\{w_{2}t_{1},w_{2}t_{2}\}$, $D_{w_{1}}=g_{11}\setminus \{x,w_{1},t\}$, $D_{w_{2}}=g_{21}\setminus \{x,w_{2},t_{1}\}$, and $D'_{w_{2}}=g_{22}\setminus \{x,w_{2},t_{2}\}$. Now let $W_{1}\cap W_{2}\neq \emptyset$. Without any loss of generality we may assume that $g_{11}=g_{21}$. Choose $t_{1}\in g_{22}\setminus \{x,w_{1}, w_{2}\}$ and set $E'_{6}=\{w_{1}w_{2}\}$, $E''_{6}=\{w_{2}w_{1},w_{2}t_{1}\}$, $D_{w_{1}}=g_{11}\setminus \{x,w_{1},w_{2}\}$, $D_{w_{2}}=g_{22}\setminus \{x,w_{1},w_{2}, t_{1}\}$ and $D'_{w_{2}}=\emptyset$. 
  \item[$\bullet$]
$E_{7}=\{ xv\vert v\in (V(\Gamma)\setminus (\{x,y,z\}\cup D))\cup \{w_{1}, w_{2}\}\}$.
\end{itemize}
\begin{emp}\label{Hhamilt1}
The graph  $\Gamma$ is Hamiltonian.
\end{emp}

Assume that $d_1\leq d_2\leq \cdots\leq d_n$ are degrees of the vertices of $\Gamma$. Now we  show that for each $i\leq \frac{n}{2}$, we have $d_{i}>i$ or  $d_{n-i}\geq n-i$. So Chv\'{a}tal's condition  \cite{chvat} implies the existence of a Hamiltonian cycle in $\Gamma$.
Clearly, $d_{\Gamma}(x)\geq n-6$. When $u\in U_{12}\setminus \{y\}$, apart from at most four choices of $v \in V(\Gamma)\setminus \{u,x,y,z\}$ the edges $\{x,z,u,v\}$ of $H$ are of color $1$. So $d_{\Gamma}(u)\geq n-8$ where $u\in U_{12}\setminus \{y\}$. Similarly, $d_{\Gamma}(u)\geq n-7$ for $u\in U_{13}\setminus \{z\}$ and also we have  $d_{\Gamma}(u)\geq n-6$ when  $u\in U_{13}\setminus (\{z\}\cup D)$. It is  straightforward to see that $d_{\Gamma}(z)\geq n-\vert A\vert -\vert A'\vert -7\geq \frac{n+1}{2}$ and $d_{\Gamma}(y)\geq n-\vert B\vert -\vert B'\vert -7\geq \frac{n+3}{2}$. For $U_{123}=\emptyset$, Chv\'{a}tal's condition implies that the  graph $\Gamma$ is Hamiltonian. Now let $U_{123}=\{w_{1}, w_{2},\ldots ,w_{m}\}\neq \emptyset$, $\vert U_{12}\setminus \{y\}\vert =l$, $\vert U_{13}\vert =k$ and suppose that $d_{\Gamma}(w_{i})\leq d_{\Gamma}(w_{i+1})$ for every $1\leq i\leq m-1$.
For $i=1,\ldots ,m$, let
 $$N_{i}=\{ \{ x,z,v,w_{i}\} \vert v\in U_{12}\setminus \{y\}\}\cup \{ \{ x,y,v,w_{i}\} \vert v\in U_{13}\}.$$
For each $1\leq i\leq m$, suppose that $n_{i}$ is the number of edges of color $1$ in $N_{i}$.
Clearly, $d_{\Gamma}(w_{i})\geq n_{i}$, for each $1\leq i\leq m$. Among all $m(k+l)$ edges in $\bigcup _{i=1}^{m}N_{i}$, there are at most $2(k+l)+2$ edges  of colors $2$ and $3$. So $\sum_{i=1}^{m}n_{i} \geq (m-2)(k+l)-2$.
If $d_{\Gamma}(w_{3})\leq \lfloor \frac{k+l}{3}\rfloor -1$, then $\sum_{i=1}^{3} n_{i}\leq \sum_{i=1}^{3} d_{\Gamma}(w_{i}) \leq k+l-3$. Therefore,
$$\sum_{i=4}^{m} n_{i}\geq   (m-2)(k+l)-2-(k+l-3)=(m-3)(k+l)+1,$$
which is impossible, since $\vert \bigcup _{i=4}^{m}N_{i}\vert =(m-3)(k+l) $.
Thus, $d_{\Gamma}(w_{3})\geq \lfloor \frac{k+l}{3}\rfloor > 9$ and consequently $d_{\Gamma}(w_{i})\geq 10 $ for $ 3\leq i\leq 6$. On the other hand, by the definitions of $E_{6}$ and $E_{7}$ we have $d_{\Gamma}(w_{1})\geq 2$ and $d_{\Gamma}(w_{2})\geq 3$. Hence,
\begin{equation}\label{ff}
d_{i}>i  \ \ \  $ for $ 3\leq i\leq 6.
\end{equation}
Since $\vert U_{123}\vert < \frac{n}{2}$ and $\vert U_{13}\vert \geq \frac{1}{2}\lfloor \frac{n}{2} \rfloor $, we have  $d_{n-i}\geq n-i$ for each $6 \leq i\leq \frac{n}{2}$. On the other hand, by (\ref{ff}), we have $d_{i}>i$ for $1\leq i\leq 6$. Now clearly Chv\'{a}tal's condition  yields the existence of a Hamiltonian cycle in $\Gamma$.
\begin{emp}\label{extend1}
Every Hamiltonian cycle in  $\Gamma$ can be extended to a monochromatic Hamiltonian Berge-cycle of color $1$ in $H$.
\end{emp}

Suppose that $v_{1}, v_{2},\ldots ,v_{n-1},v_{n}=x$ is the vertices of a Hamiltonian cycle $C$ in $\Gamma$. Without any loss of generality, we may assume that $v_1\neq w_1$. Now for $i=1,2,\ldots,n$, we define the edges $f_{i}\in E(H)$ of color $1$ in the same order their subscripts appear so that $\{v_{i},v_{i+1}\}\subseteq  f_{i}$ and $f_{1}, f_{2},\ldots ,f_{n}$  make a Hamiltonian Berge-cycle with the core sequence $v_{1}, v_{2}, \ldots ,v_{n}$. First let $i=1,2,\ldots,n-2$. Set $f_{i}=\{x,z,v_i ,v_{i+1}\}$ for  $v_{i}v_{i+1}\in E_{1}\cup E_{4}$ and  $f_{i}=\{x,y,v_i ,v_{i+1}\}$ for $v_{i}v_{i+1}\in E_{2}\cup E_{3}$. If $v_{i}v_{i+1}\in E_{5}$, then set  $f_{i}=\{x,v_{i},v_{i+1}, u\}$ of color $1$ so that $u\in U_{13}\setminus \{z,v_{i-1},v_{i},v_{i+1},v_{i+2}, v_{1},v_{n-1}\}$. Such an edge exists since  $\vert U_{13}\vert \geq \frac{1}{2} \lfloor \frac{n}{2}\rfloor\geq 20$ and for a fixed vertex $v\in U_{13}\setminus \{z\}$ there are at least  $q=\frac{1}{2} \lfloor \frac{n}{2}\rfloor -6 > 14$ vertices, say $\{u_1 ,u_2 ,\ldots ,u_{q}\}$ in $U_{13}\setminus \{z, v\}$, where every edge $\{x,y,v,u_{j}\}$ is of color $1$. If $v_{i}v_{i+1}\in E_{6}$, then  by the definition of $E_{6}$, there is an appropriate edge  $f_{i}\in W_{1}\cup W_{2}$  containing $v_{i}$ and $v_{i+1}$.
 Now let $i=n-1 $. It is easy to see that  $\{ v_{n-1},x\}$  has been used in at most two of the edges  $f_{i}$s for $1\leq i\leq n-2$. On the other hand,  $1\in c_2^{*}(v_{n-1}x)$. Thus  we can choose an appropriate edge $f_{n-1}$. Finally let $i=n$. One can see that  $\{x, v_{1}\}$  has been used in at most two of the edges  $f_{i}$s for $1\leq i\leq n-1$  and since  $1\in c_2^{*}(xv_{1})$, then  there is an appropriate edge $f_{n}$.

\bigskip
\noindent \textbf{Case 2. }$\vert U_{23}\vert =1$.\\

Since $\vert U_{23}\vert =1$, we have $1\leq \vert U_{12}\vert \leq 2$. Assume that $U_{23}=\{u_{23}\}$, 
$U_{12}=\{y, u_{12}\}$ for $|U_{12}|=2$ and $U_{12}=\{y\}$ for $|U_{12}|=1$. If $\frac{n-2}{2} \leq\vert U_{123}\vert \leq \frac{n-1}{2}$ and
 $U_{123}\subseteq B_{1}$, then for each $w\in U_{123}$ and each $v\in V(G)$ we have $2\in c_2^{*}(wv)$. On the other hand, $2\in c_2^{*}(x u_{23})$ and $2\in c_2^{*}(xy)$ and so  Chv\'{a}tal's condition implies that the subgraph induced by all edges $e$ with $2\in c_2^{*}(e)$   contains a Hamiltonian cycle. By Lemma \ref{bb} for $r=4$ and $t=2$, the proof is completed.\\
Now fix a  vertex $w\in U_{123}\setminus B_{1}$ when $\frac{n-2}{2} \leq\vert U_{123}\vert \leq \frac{n-1}{2}$. Let $U_{123}$ be partitioned into $A,B$, where $\vert B\vert \leq \vert A| \leq \vert B\vert +1$. Consider a graph $\Gamma$ with the vertex set $V(\Gamma)=V(H)$, and the edge set
$E(\Gamma)=\bigcup _{i=1}^{8}E_{i}$, where $E_{i}$s are defined as follows:
\begin{itemize}
\item[$\bullet$]
$E_{1}=\emptyset$ when $U_{12}=\{y\}$ and $E_{1}=\{ u_{12}v\vert v\neq y, c(\{ x,z,u_{12},v\})=1\}$, otherwise.
\item[$\bullet$]
$E_{2}=\{ uv\vert u\in U_{13}\setminus \{z\}, c(\{ x,y,u,v\})=1\}$.
\item[$\bullet$]
$E_{3}=\{ zv\vert v\in V(\Gamma)\setminus A, c(\{ x,y,z,v\})=1\}$.
\item[$\bullet$]
$E_{4}=\{ yv\vert v\in A, c(\{ x,y,z,v\})=1\}$.
\item[$\bullet$]
$E_{5}=\{ yv\vert v\in U_{13}\setminus \{z\}\}$.
\item[$\bullet$]
Assume that $U_{123}=\{w_{1}, w_{2},\ldots ,w_{m}\}$ and $d_{\Gamma'}(w_{1})\leq d_{\Gamma'}(w_{2})\leq \cdots \leq d_{{\Gamma'}}(w_{m})$, where $\Gamma'$ is the graph induced by
 $\bigcup _{i=1}^{5} E_{i}$. Assume that $e_{w_{1}v}=\{x,z,w_{1},v\}$ (resp. $e_{u_{23}v}=\{x,z,u_{23},v\}$) when $w_{1}v \in E_{1}\cup E_{4}$ (resp. $u_{23}v \in E_{1}$) and $e_{w_{1}v}=\{x,y,w_{1},v\}$ (resp. $e_{u_{23}v}=\{x,y,u_{23},v\}$) when $w_{1}v \in E_{2}\cup E_{3}$ (resp. $u_{23}v \in E_{2}\cup E_{3}$).  By Lemma \ref{a} and the fact $1\in c_2^{*}(xw_{1})$, there are $r=$ max $\{3-d_{{\Gamma'}}(w_{1}),0\}$ and $l=$ max $\{2-d_{{\Gamma'}}(u_{23}),0\}$  edges $W=\{h_{1}, \ldots ,h_{r}\} \subseteq E(H)\setminus \{e_{w_{1}v} |w_{1}v\in \bigcup_{i=1}^{4} E_{i}\}$  and $U=\{g_{1}, \ldots ,g_{l}\} \subseteq E(H)\setminus \{e_{u_{23}v} |u_{23}v\in \bigcup_{i=1}^{3} E_{i}\}$
 of color $1$ containing $\{x,w_{1}\}$ and $u_{23}$,  respectively. We consider three cases:
 \begin{itemize}
 \item[i.] 
 $r\leq 1$. Set $E'_{6}=D_{w_{1}}=D'_{w_{1}}=\emptyset$.  If $l=0$, then set $ E''_{6}=D_{u_{23}}=D'_{u_{23}}=\emptyset$. If $l=1$, then $U' = \{g_{1}, e_{u_{23}v}\}$ where $u_{23}v\in  \bigcup_{i=1}^{3} E_{i}$. Set $D'_{u_{23}}=\emptyset$, $E''_{6}=\{u_{23}t_{1}\}$, $D_{u_{23}}=g_{1}\setminus \{x,u_{23},t_{1}\}$ so that $t_{1}\in g_{1}\setminus \{x,u_{23},v\}$.  Now  let $l=2$. Thus $U= \{g_{1}, g_{2}\}$. Set  $E''_{6}=\{u_{23}t_{1},u_{23}t_{2}\}$, $D_{u_{23}}=g_{1}\setminus \{x,u_{23},t_{1}\}$ and $D'_{u_{23}}=g_{2}\setminus \{x,u_{23},t_{2}\}$ so that $t_{1}\in g_{1}\setminus \{x,u_{23}\}$  and $t_{2}\in g_{2}\setminus \{x,u_{23}, t_{1}\}$  and $t_{1}$ and $t_{2}$ have maximum repetitions in $g_{1}$ and $g_{2}$.
 \item[ii.]
 $r=2$. Set $D'_{w_{1}}=\emptyset$. Let $W'=\{h_{1}, h_{2},e_{w_{1}v}\}$ where $w_{1}v\in \bigcup_{i=1}^{4} E_{i}$. If $l=0$, then set $E''_{6}=D_{u_{23}}=D'_{u_{23}}=\emptyset$, $E'_{6}=\{w_{1}u_{1}\}$, $D_{w_{1}}=h_{1}\setminus \{x,w_{1},u_{1}\}$ where $u_{1}\in h_{1}\setminus \{x,w_{1}, v\}$. If $l=1$, then $U' = \{g_{1}, e_{u_{23}u}\}$ where $u_{23}u\in  \bigcup_{i=1}^{3} E_{i}$. We may assume that $\vert h_{1}\cap \{u_{23}\}\vert \geq \vert h_{2}\cap \{u_{23}\}\vert $. Set $D'_{u_{23}}=\emptyset$, $E'_{6}=\{w_{1}u_{1}\}$, $E''_{6}=\{u_{23}t_{1}\}$, $D_{w_{1}}=h_{2}\setminus \{x,w_{1},u_{1}\}$, $D_{u_{23}}=g_{1}\setminus \{x,u_{23},t_{1}\}$ so that  $u_{1}\in h_{2}\setminus \{x,w_{1},v\}$, $t_{1}\in g_{1}\setminus \{x,u_{23},u\}$.  Now let $l=2$. Then we have $U=\{g_{1},g_{2}\}$. If $W'\cap U=\emptyset$, then set  $E'_{6}=\{w_{1}u_{1}\}$, $E''_{6}=\{u_{23}t_{1}, u_{23}t_{2}\}$, $D_{w_{1}}=h_{1}\setminus \{x,w_{1},u_{1}\}$, $D_{u_{23}}=g_{1}\setminus \{x,u_{23},t_{1}\}$ and  $D'_{u_{23}}=g_{2}\setminus \{x,u_{23},t_{2}\}$ so that $u_{1}\in h_{1}\setminus \{x,w_{1},v\}$, $t_{1}\in g_{1}\setminus \{x,u_{23}\}$,  $t_{2}\in g_{2}\setminus \{x,u_{23}, t_{1}\}$ and $t_{1}$ and $t_{2}$ have maximum  repetitions in $g_{1}$ and $g_{2}$. Otherwise, we may assume that $h_{1}=g_{1}$. Set  $D'_{u_{23}}=\emptyset$, $E'_{6}=\{w_{1}u_{23}\}$, $E''_{6}=\{u_{23}w_{1},u_{23}t_{1}\}$, $D_{w_{1}}=h_{1}\setminus \{x,w_{1},u_{23}\}$, $D_{u_{23}}=g_{2}\setminus \{x,u_{23},t_{1}\}$ where $t_{1}\in g_{2}\setminus \{x,u_{23}, w_{1}\}$ and $t_{1}$ has maximum  repetition in $g_{1}$ and $g_{2}$.
 \item[iii.]
 $r=3$.  If $l=0$, then set $E''_{6}=D_{u_{23}}=D'_{u_{23}}=\emptyset$, $E'_{6}=\{w_{1}u_{1}, w_{1}u_{2}\}$, $D_{w_{1}}=h_{1}\setminus \{x,w_{1},u_{1}\}$, $D'_{w_{1}}=h_{2}\setminus \{x,w_{1},u_{2}\}$ where $u_{1}\in h_{1}\setminus \{x,w_{1}\}$, $u_{2}\in h_{2}\setminus \{x,w_{1},u_{1}\}$ and $u_{1}$ and $u_{2}$ have maximum repetitions in $h_{1}$ and $h_{2}$. If $l=1$, then $U'= \{g_{1}, e_{u_{23}v}\}$ where $u_{23}v\in  \bigcup_{i=1}^{3} E_{i}$. We may assume that $g_{1}\notin \{h_{2},h_{3}\}$. Now set $D'_{u_{23}}=\emptyset$, $E'_{6}=\{w_{1}u_{1}, w_{1}u_{2}\}$, $E''_{6}=\{u_{23}t_{1}\}$, $D_{w_{1}}=h_{2}\setminus \{x,w_{1},u_{1}\}$, $D'_{w_{1}}=h_{3}\setminus \{x,w_{1},u_{2}\}$, $D_{u_{23}}=g_{1}\setminus \{x,u_{23},t_{1}\}$  where $u_{1}\in h_{2}\setminus \{x,w_{1}\}$, $u_{2}\in h_{3}\setminus \{x,w_{1},u_{1}\}$, $t_{1}\in g_{1}\setminus \{x,u_{23}, v\}$ and $u_{1}$ and $u_{2}$ have maximum repetitions in $h_{2}$ and $h_{3}$. Now let $l=2$.  If $W\cap U=\emptyset$, then set  $E'_{6}=\{w_{1}u_{1}, w_{1}u_{2}\}$, $E''_{6}=\{u_{23}t_{1},u_{23}t_{2}\}$, $D_{w_{1}}=h_{1}\setminus \{x,w_{1},u_{1}\}$, $D'_{w_{1}}=h_{2}\setminus \{x,w_{1},u_{2}\}$, $D_{u_{23}}=g_{1}\setminus \{x,u_{23},t_{1}\}$, $D'_{u_{23}}=g_{2}\setminus \{x,u_{23},t_{2}\}$ so that   $u_{1}\in h_{1}\setminus \{x,w_{1}\}$, $u_{2}\in h_{2}\setminus \{x,w_{1}, u_{1}\}$ are vertices with maximum repetitions in $h_{1}$ and $h_{2}$ and   $t_{1}\in g_{1}\setminus \{x,u_{23}\}$  and $t_{2}\in g_{2}\setminus \{x,u_{23}, t_{1}\}$ are vertices with maximum repetitions in $g_{1}$ and $g_{2}$. Otherwise, we may assume that $\vert h_{i}\cap \{u_{23}\}\vert \geq \vert h_{j}\cap \{u_{23}\}\vert $ for $i<j$, $h_{1}=g_{1}$ and $h_{3}\notin U$. Set $D'_{u_{23}}=\emptyset$, $E'_{6}=\{w_{1}u_{23}, w_{1}u_{1}\}$, $E''_{6}=\{u_{23}w_{1},u_{23}t_{1}\}$, $D_{w_{1}}=h_{1}\setminus \{x,w_{1},u_{23}\}$, $D'_{w_{1}}=h_{3}\setminus \{x,w_{1}, u_{1}\}$, $D_{u_{23}}=g_{2}\setminus \{x,u_{23},w_{1}, t_{1}\}$ so that $u_{1}\in h_{3}\setminus \{x,w_{1}\}$, $t_{1}\in g_{2}\setminus \{x,u_{23}, w_{1}\}$ and $t_{1}$ has maximum repetition in $g_{1}$ and $g_{2}$. 
 \end{itemize}
 In all cases set $E_{6}=E'_{6}\cup E''_{6}$ and $D=D_{w_{1}}\cup D'_{w_{1}}\cup D_{u_{23}}\cup D'_{u_{23}}$.\\
 \item[$\bullet$]
$E_{7}=\emptyset $ if $\vert U_{123}\vert \leq \frac{n-3 }{2}$ and $E_{7}=\{wv \vert \  v\in V(\Gamma)\setminus \{x,w\},  1\in c_2^{*}(vw)\}$, otherwise. It is easy to see that $|E_{7}|\geq \frac{n}{2}$ when  $\frac{n-2}{2} \leq\vert U_{123}\vert \leq \frac{n-1}{2}$.
 \item[$\bullet$]
$E_{8}=\{ xv\vert v\in (V(\Gamma)\setminus (\{x,y,z,u_{23},w\}\cup D)) \cup \{w_{1}\}\}$.
\end{itemize}

\begin{emp}\label{Hhamilt2}
The graph  $\Gamma$ is Hamiltonian.
\end{emp}
Assume that $d_1\leq d_2\leq \cdots\leq d_n$ are degrees of the vertices of $\Gamma$.  Now we  show that for each $i\leq \frac{n}{2}$, we have $d_{i}>i$ or  $d_{n-i}\geq n-i$. So Chv\'{a}tal's condition  implies the existence of a Hamiltonian cycle in $\Gamma$.  According to the above discussions  $d_{\Gamma}(x)\geq n-11$. For every $u\in U_{13}\setminus \{z\}$, with at most four choices of $v \in V(\Gamma)\setminus \{x,y,u\}$ excluded the edges $\{x,y,u,v\}$ of $H$ are of color $1$. So
 $d_{\Gamma}(u)\geq n-7$, where $u\in U_{13}\setminus \{z\}$ and also we have  $d_{\Gamma}(u)\geq n-6$ when  $u\in U_{13}\setminus ( \{z\}\cup D) $. Similarly, $d_{\Gamma}(u_{12})\geq n-8$ when $U_{12}=\{y, u_{12}\}$. It is straightforward to see that $d_{\Gamma}(u_{23})\geq 2$,   $d_{\Gamma}(z)\geq n-\vert A\vert -7 \geq \frac{n+5}{2}$ and $d_{\Gamma}(y)\geq n-\vert B\vert -9 \geq \frac{n+3}{2}$. If $U_{123}=\emptyset$, then one can easily see that Chv\'{a}tal's condition implies that the  graph $\Gamma$ is Hamiltonian.\\
  Now let  $U_{123}=\{w_{1}, w_{2},\ldots ,w_{m}\}\neq \emptyset$ with $d_{\Gamma}(w_1)\leq d_{\Gamma}(w_2)\leq \cdots\leq d_{\Gamma}(w_m)$ and $\vert U_{13}\vert =k$.
We show that $d_{\Gamma}(w_{m})\geq  \frac{n}{2}$ when $\frac{n-2}{2} \leq\vert U_{123}\vert \leq \frac{n-1}{2}$. If $w\in \bigcup _{i=2}^{4}B_{i}$, then one can easily see that $\vert E_{7}\vert=\vert \{v\vert \  v\in V(\Gamma)\setminus \{x,w\}, 1\in c_2^{*}(wv)\}\vert \geq \frac{n}{2}$. So $d_{\Gamma}(w_{m})\geq d_{\Gamma}(w) \geq\frac{n}{2}$. Now let $w\in V(G)\setminus  \bigcup _{i=1}^{4}B_{i}$. From the definition of $x$, we conclude that  $\frac{n-2}{2} \leq\vert U_{123}(w)\vert \leq \frac{n-1 }{2}$ and $U_{13}(w)$ and $U_{12}(w)$ are non-empty. Hence again $\vert E_{7}\vert=\vert \{v\vert \ v\in V(\Gamma)\setminus \{x,w\}, 1\in c_2^{*}(wv)\}\vert \geq \frac{n}{2}$ and so we have $d_{\Gamma}(w_{m})\geq d_{\Gamma}(w) \geq\frac{n}{2}$ when $\frac{n-2}{2} \leq\vert U_{123}\vert \leq \frac{n-1}{2}$.

 Now we claim that
\begin{equation}\label{pp}
d_{i}>i \ \ \ \ \ $for$ \ \ \ 1\leq i\leq 6.
\end{equation}
First let $U_{12}=\{y, u_{12}\}$ and let $T=\{ \{x,y,u_{23},v\}, \{x,u_{12},u_{23},v\} \vert v\in U_{13}\}.$
Let $S$ be the set of all vertices $v\in U_{13}$ for which there is an edge of color 1 or 3 containing $v$ in $T$. Clearly, $|S|\leq 6$. Therefore, for each $v\in U_{13}\setminus S$, since $2\notin c_2^{*}(xv)$ apart from two edges in $T$ all edges of $H$ containing $x$ and $v$ are of color $1$ or $3$. On the other hand, since $3\notin c_2^{*}(xy)$ at most two edges $\{x,y,v,w_{i}\}$ are of color $3$ where $i\in \{1,2,\ldots ,m\}$. So
for each $i=1,\ldots ,m$ at least $k-8 >30$ edges in $\{ \{ x,y,v,w_{i}\} \vert v\in U_{13}\setminus S\}$ are of color 1. Hence, $d_{\Gamma}(w_{1}) >30$ and so $d_{i}>i$ for each  $1\leq i\leq 6$.

Now let $U_{12}=\{y\}$ and $T=\{ \{x,y,u_{23},v\}\vert v\in U_{13}\}.$ At least $k-4$ edges in $T$ are of color $2$ in $H$. Now
for $i=1,\ldots ,m$, consider $N_{i}=\{ \{ x,y,v,w_{i}\} \vert v\in U_{13}\}.$ For every $1\leq i\leq m$, suppose that $n_{i}$ is the number of edges of color $1$ in $N_{i}$.
Clearly, $d_{\Gamma}(w_{i}) \geq n_{i}$. Let $S_{T}\subseteq U_{13}$ be the set of all vertices $v$ that lies on an edge of $T$ of color $2$. Clearly, $|S_{T}|\geq k-4$. Since $2\notin c_2^{*}(xv)$ for each $v\in U_{13}$, there are at most $|S_{T}|$ (resp. $2(k-|S_{T}|)$) edges of color $2$ in  $\bigcup _{i=1}^{m}N_{i}$ each containing a vertex in $S_{T}$ (resp. $U_{13}\setminus S_{T}$). Therefore,  among all $mk$ edges in $\bigcup _{i=1}^{m}N_{i}$ there are at most $k+6$ edges  of colors $2$ and $3$. So $\sum_{i=1}^{m}n_{i}\geq (m-1)k-6.$
If $d_{\Gamma}(w_{2})\leq \lfloor \frac{k-7}{2} \rfloor$,  then $\sum_{i=1}^{2} n_{i}\leq \sum_{i=1}^{2} d_{H}(w_{i})\leq k-7$. Therefore, 
$$\sum_{i=3}^{m}n_{i} \geq (m-1)k-6-(k-7)= (m-2)k+1,$$
which is impossible, since $\vert \bigcup _{i=3}^{m}N_{i}\vert =(m-2)k$.
Thus, $d_{\Gamma}(w_{2})> \lfloor \frac{k-7}{2} \rfloor\geq15$ and consequently $d_{\Gamma}(w_{i})\geq 16$   for  $2\leq i\leq 6$. On the other hand,  according to the definitions of $E_{6}$ and $E_{8}$, we have $d_{u_{23}}\geq 2$ and $d_{w_{1}}\geq 3$. Therefore, 
$d_{i}>i$ for each  $1\leq i\leq 6$.

Based on the previous discussions, since $\vert U_{123}\vert \leq \frac{n-1}{2}$ and $|U_{13}| \geq \lfloor \frac{n}{2} \rfloor -3$, we have  $d_{n-i}\geq n-i$ for each $6 \leq i\leq \frac{n}{2}$. On the other hand by (\ref{pp}), we have $d_{i}>i$ for $1\leq i\leq 6$. Now, Chv\'{a}tal's condition implies the existence of a Hamiltonian cycle in $H$.

\begin{emp}\label{extend2}
There is a Hamiltonian Berge-cycle of color $1$ in $H$.
\end{emp}

We show that every Hamiltonian cycle in $\Gamma$ can be extended to a monochromatic Hamiltonian Berge-cycle in $H$.
Suppose that $v_{1}, v_{2},\ldots ,v_{n-1},v_{n}=x$ is the vertices of a Hamiltonian cycle $C$ in $\Gamma$. 
%Without any loss of generality, we may assume that $v_1\neq w_1$.
 Now for each $i=1,2,\ldots ,n$, we define an edge $f_{i}\in E(H)$ of color $1$ in the same order  their subscripts appear so that $\{v_{i},v_{i+1}\}\subseteq  f_{i}$ and $f_{1}, f_{2},\ldots ,f_{n}$  make a Hamiltonian Berge-cycle with the core sequence $v_{1}, v_{2}, \ldots ,v_{n}$.
First let $i\in [n]\setminus (\{n-1,n\}\cup \{i\vert v_{i}v_{i+1}\in E_{7}\})$, where $[n]=\{1,2,\ldots ,n\}$. Set $f_{i}=\{x,z,v_i ,v_{i+1}\}$ for  $v_{i}v_{i+1}\in  E_{1}\cup E_{4}$ and  $f_{i}=\{x,y,v_i ,v_{i+1}\}$ for $v_{i}v_{i+1}\in E_{2}\cup E_{3}$. Now let $v_{i}v_{i+1}\in E_{5}$. Set $f_{i}=\{x,v_{i},v_{i+1},u\}$ of color $1$, where $u\in U_{13}\setminus \{z,v_{i-1},v_{i},v_{i+1},v_{i+2}, v_{1},v_{n-1}\}$. Such an edge exists since  $\vert U_{13}\vert \geq \lfloor \frac{n}{2} \rfloor -3$ and for a fixed vertex $v\in U_{13}\setminus \{z\}$ there are at least  $q=\lfloor \frac{n}{2} \rfloor -9 >30$ vertices, say $\{u_1 ,u_2 ,\ldots ,u_{q}\}$ in $U_{13}\setminus \{z, v\}$, where every edge $\{x,y,v,u_{j}\}$ is of color $1$.  If $v_{i}v_{i+1}\in E_{6}$, then  by the definition of $E_{6}$, there is an appropriate edge  $f_{i}\in W\cup U$  containing $v_{i}$ and $v_{i+1}$.

Now let $L_{uv}\subset E(H)\setminus \{f_{i}\vert \  i\in [n]\setminus (\{n-1,n\}\cup \{i\vert \ v_{i}v_{i+1}\in E_{7}\})\}$ be the set of all edges of color $1$ containing $u$ and $v$. Note that  
$1\in c_2^{*}(v_{n-1}x)$ and $1\in c_2^{*}(xv_{1})$. By the definitions of $E_{6}$ and $E_{7}$, it is easy to see that Hall's theorem implies that we can choose appropriate edges $f_{n-1}\in L_{v_{n-1}x}$, $f_{n}\in L_{xv_{1}}$ and $f_{i}\in L_{v_{i}v{i+1}}$ for each $i$ with $ v_{i}v_{i+1}\in E_{7}$.

\epf
%%%%%%%%%%%%%%%%%%%%%%%%%%%%%%%%%%%%%%%%%%%%%%%%%%%%%%%%%555
%%%%%%%%%%%%%%%%%%%%%%%%%%%%%%%%%%%%%%%%%%%%%%%%%%%%%%%555


\footnotesize
\begin{thebibliography}{99}

\bibitem{Berge} C. Berge, Graphs and Hypergraphs, North Holland, Amsterdam and London, 1973.

\bibitem{DGS}  P. Dorbec, S. Gravier and G.N. S\'{a}rk\"{o}zy, Monochromatic Hamiltonian $t$-tight Berge-cycles in hypergraphs,
{\em J. Graph Theory} {\bf 59} (2008), 34--44.

\bibitem{gyarf1} A. Gy\'{a}rf\'{a}s, J. Lehel, G.N. S\'{a}rk\"{o}zy and R.H. Schelp, Monochromatic Hamiltonian Berge-cycles in colored complete uniform hypergraphs, {\em J. Combin. Theory Ser. B. } {\bf 98} (2008), 342--358.

\bibitem{gyarf2} A. Gy\'{a}rf\'{a}s, G.N. S\'{a}rk\"{o}zy and E. Szemer\'{e}di, Long monochromatic Berge-cycles in colored $4$-uniform hypergraphs, {\em Graphs Combin.} {\bf 26} (2010), 71--76.

\bibitem{GSS2} A. Gy\'{a}rf\'{a}s, G.N. S\'{a}rk\"{o}zy and E. Szemer\'{e}di, Monochromatic
Hamiltonian $3$-tight Berge cycles in $2$-colored $4$-uniform hypergraphs, {\em J. Graph Theory } {\bf 63} (2010), 288--299.

\bibitem{GSS}  A. Gy\'{a}rf\'{a}s, G.N. S\'{a}rk\"{o}zy and E. Szemer\'{e}di, Monochromatic matchings in the shadow graph of almost
complete hypergraphs, {\em Ann. Combin.}  {\bf 14} (2010), 245--249.


\bibitem{hax1} P. Haxell, T. Luczak, Y. Peng, V. R\"{o}dl, A. Ruci\'{n}ski, M. Simonovits and  J. Skokan, The Ramsey number for hypergraph cycles I, {\it J. Combin. Theory Ser. A} {\bf 113} (2006), 67--83.

\bibitem{hax2}  P. Haxel, T. Luczak, Y. Peng, V. R\"{o}dl, A. Ruci\'{n}ski and J. Skokan, The Ramsey number for $3$-uniform tight hypergraph cycles, {\em Combin. Probab. Comput.} {\bf 18} (2009), 165--203.

\bibitem{chvat}  L. Lov\'{a}sz,  Combinatorial Problems and Exercises, 2nd edn. Akad\'{e}miai Kiad\'{o}, North Holland (1976).

\bibitem{OS} G.R. Omidi and M. Shahsiah, Ramsey numbers of 3-uniform
loose paths and loose cycles, {\it J. Combin.
 Theory Ser. A}, {\bf 121} (2014), 64--73.


\bibitem{Ramsey}  F.P. Ramsey,  On a problem of formal logic, {\em Proc. London Math. Soc. Ser. 2} {\bf 30} (1930), 264--286.

\bibitem{RRS}  V. R\"{o}dl, A. Ruci\'{n}ski and E. Szemer\'{e}di, A dirac-type theorem for 3-uniform hypergraphs, {\em Comb. Probab.
Comput.} {\bf 15} (2006), 229--251.

\end{thebibliography}
\end{document}